\documentclass[12pt]{article}
\usepackage[cp1251]{inputenc}
\usepackage[russian]{babel}
\begin{document}
\author{Белошапка В.К.}

\date{25.02.2021}

\title{Об исключительных квадриках}

\maketitle

\begin{abstract}
 Доказано, что градуированные алгебры Ли инфинитезимальных голоморфных автоморфизмов невырожденной квадрики коразмерности $k$ не имеют ненулевых компонент веса больше $2k$ . Доказано также, что при  $k \leq 3$ нет градуированных компонент веса больше чем 2.
Сформулирован ряд вопросов.
 \end{abstract}

\footnote{
Механико-математический факультет Московского университета им.Ломоносова,
Воробьевы горы, 119992 Москва, Россия, vkb@strogino.ru }

{\bf  Введение}

\vspace{3ex}

Недавно было обнаружено \cite{BM20}, что доказательство одного из утверждений моей работы 1990-го года \cite{VB90} содержит ошибку. Речь идет о теореме на стр.19, из которой, в частности, следует, что градуированная алгебра Ли ${\rm aut} \, \mathcal{Q}$ инфинитезимальных голоморфных автоморфизмов невырожденной квадрики $\mathcal{Q}$ произвольной коразмерности  состоит не более чем из пяти весовых компонент
$$ {\rm aut} \, \mathcal{Q} = g_{-2}+g_{-1}+g_{0}+g_{1}+g_{2}.$$
В \cite{BM20} к этому утверждению был предъявлен контрпример квадрики коразмерности пять в ${\bf C}^9$, а среди  контрпримеров работы \cite{GM} имеется квадрика коразмерности четыре в ${\bf C}^{10}$.  Такие квадрики, т.е. невырожденные квадрики чья алгебра содержит поля веса больше чем два, мы называем {\it исключительными}. Все имеющиеся примеры исключительных квадрик весьма сложны. Природа этих редких примеров по сей день представляется неясной.

Работа построена следующим образом.   В разделе 1 анализируем условия принадлежности поля алгебре автоморфизмов в контексте теоремы Эренпрайса-Паламодова.Это та техника работы с  системами линейных дифференциальных уравнений с постоянными коэффициентами,  анализ которых позволил  сформулировать критерий конечномерности алгебры Ли инфинитезимальных голоморфных автоморфизмов квадрики  \cite{VB88}. В теореме 4 приводится критерий исключительности квадрики.  В разделе 2 мы применяем полученый критерий к анализу квадрик $CR$-типов $(3,3)$ и $(3,4)$. Использование классификаций, построенных в работах Н.Палинчак \cite{NP} и С.Анисовой \cite{SA} позволяет показать, что среди квадрик этих типов исключительных нет. В разделе 3 обсуждается
ситуация с квадриками малых коразмерностей. Основной результат -- теорема 19 -- утверждение, что исключительных квадрик в коразмерности три -- не существует. В разделе 4 происходит возвращение в контекст раздела 1.  Квадрике сопоставляется некий подмодуль свободного модуля над кольцом полиномов (характеристический подмодуль),  в его терминах приводится еще один критерий исключительности квадрики, а также, в этих терминах, доказывается верхняя оценка на веса и степени полей из алгебры автоморфизмов квадрики через ее коразмерность.

\vspace{3ex}

{\bf 1. Критерий исключительности}

\vspace{3ex}

Пусть $M_{\xi}$ -- росток гладкого вещественного порождающего подмногообразия  комплексного пространства $\mathbf{C}^N$ ($N=n+k$) $CR$-размерности $n>0$ и коразмерности $k>0$. Пусть
$$z=(z_1,\dots, z_n), \quad w=(w_1, \dots, w_k),  \quad w_j=u_j+i \, v_j, \quad j=1, \dots,k, \; - $$
координаты в $\mathbf{C}^N$. Назначим веса переменным:
$$[z]=[\bar{z}]=1, \quad [w]=[\bar{w}]=[u]=2.$$
Это позволяет разлагать степенные ряды на весовые компоненты. А если считать, что дифференцирование по $z$ и $\bar{z}$ имеет вес $(-1)$, а дифференцирование по $w$ и $\bar{w}$ -- вес $(-2)$, то
градуированной становится и алгебра Ли формальных векторных полей, и все ее подалгебры.
После простого квадратичного преобразования уравнение $M_{\xi}$ можно записать в виде
$$M_{\xi} = \{v_j=<z,\bar{z}>_j + o(2), \quad j=1, \dots,k\}=\{v= <z,\bar{z}>+o(2)\},$$
где $<z,\bar{z}>_j$ -- эрмитовы формы на $\mathbf{C}^n$, а $o(m)$ -- это функции, тейлоровское разложение которых в нуле не содержит членов веса $m$ и ниже. Основной объект нашего внимания -- это касательная квадрика
$$ \mathcal{Q}=\{v= <z,\bar{z}>\}.$$
Пусть ${\rm aut} \, M_{\xi}$ -- это алгебра Ли ростков векторных полей в $\mathbf{C}^N$, касательных к $M_{\xi}$,   вида
\begin{eqnarray} \label{f}
X=2 \, {\rm Re} \, \left( f(z,w) \,\frac{\partial}{\partial z}
+ g(z,w) \,\frac{\partial}{\partial w}  \right) =(f(z,w),g(z,w)),
\end{eqnarray}
где $f$ и $g$ -- ростки вектор-функций, голоморфных в $\xi$. Это поля, которые порождают локальные 1-параметрические группы голоморфных преобразований $M_{\xi}$.
Модельная поверхность $\mathcal{Q}$ взвешенно однородна (задается взвешенно однородными уравнениями) и голоморфно однородна (группа голоморфных автоморфизмов действует на $\mathcal{Q}$ транзитивно). Поэтому вместо ${\rm aut} \, \mathcal{Q}_{\xi}$ мы можем писать ${\rm aut} \, \mathcal{Q}$ (нет зависимости алгебры от точки).

В соответствии с нашей градуировкой
$${\rm aut} \, \mathcal{Q}= g_{-2}+g_{-1}+g_{0}+g_{1}+g_{2}+\dots \; .$$
Нетрудно проверить, что принадлежность векторного поля $X=\sum X_j$ алгебре Ли ${\rm aut} \, \mathcal{Q}$ равносильна принадлежности этой алгебре каждой его весовой компоненты $X_j$. В такой ситуации условие конечномерности алгебры равносильно условию ее конечной градуированности.  В \cite{VB88} было доказано, что критерием конечномерности ${\rm aut} \, \mathcal{Q}$ является следующая пара условий:
(1) отсутствие ядра, т.е. если $<e,\bar{z}>=0$ для всех $z$, то $e=0$;  (2) координаты формы  $<z,\bar{z}>$ линейно независимы (условие {невырожденности}).

Записывая условие принадлежности поля (\ref{f}) алгебре ${\rm aut} \, \mathcal{Q}$, получаем
\begin{eqnarray} \label{88}
 2 \, {\rm Re} \, \left( i \, g (z,u+i\,<z,\bar{z}>) + 2 \, <f(z,u+i\,<z,\bar{z}>) ,\bar{z}>  \right)=0.
\end{eqnarray}

Положим $ \Delta(\varphi(u))=\partial_u \varphi(u)(< z, \bar{z}>)$.
Несложный анализ (см. \cite{VB88}) соотношения (\ref{88})  позволяет установить следующее

\vspace{3ex}

{\bf Утверждение 1:} Пара $(f,g)$ удовлетворяют соотношению (\ref{88}) тогда и только тогда, когда
$$f=a(w)+C(w)\,z +A(w)(z,z), \quad g=b(w)+2\,i\,< z, \bar{a}(w) >.$$
(зависимость от $z$ указана явно, т.е. степени ноль, один и два, а зависимость от $w$ -- аналитическая в окрестности начала),
причем $a,A,C,b$ удовлетворяют следующим соотношениям, которые распадаются на две системы уравнений
\begin{eqnarray}\label{88a}
  <A(u)(z,z),\bar{z} >=2 \,i \, < z, \Delta \bar{a}(u) >\\
 \nonumber < z, \Delta^2 \bar{a}(u) >=0.
\end{eqnarray}
 \begin{eqnarray}\label{88b}
  {\rm Im} \, b(u)=0,\\
 \nonumber  \Delta\,b(u)=2\, {\rm Re}\, < C(u) \, z, \bar{z} >,\\
  \nonumber{\rm Im}\, < \Delta^2 \,C(u) \, z, \bar{z} >=0,\\
  \nonumber \Delta^3 \,b(u)=0.
 \end{eqnarray}

Можно написать, что
$$C(u)\,z=\sum \, C_j(u)\, z_j,   \quad A(u)(z,z)=\sum \, A_{kj}(u) \, z_k \, z_j,  \quad A_{kj}=A_{jk}.$$
Тогда каждая из систем уравнений может быть записана как система линеных уравнений с постоянными коэффицентами.
Первая на  набор вектор-функций $((a(u),A_{kj}(u))$, а вторая -- $(b(u),C_j(u))$.  Нетрудно проверить, что в случае, если эрмитова вектор-форма -- невырождена, то ни одна из этих систем не имеет экспоненциальных решений с ненулевым показателем. Применяя к каждой из этих систем теорему об экспоненциальном представлении Эренпрайса-Паламодова \cite{P}, заключаем, что пространство решений каждой из них -- это конечномерное подпространство пространства полиномов от $u$.  Ясно, что конечномерность подразумевает равномерную
оценку степени решений.

При этом следует отметить, что любое векторное поле $X$ можно разложить в сумму четной -- $X^0$ и нечетной -- $X^1$ компонент.  Где --
$$X^0=X_{-2}+X_0+X_2+\dots,  \quad X^1=X_{-1}+X_1+X_3+\dots$$
 в нашей градуировке. При этом
можно заметить, что если $X \in {\rm aut} \, \mathcal{Q}$, то
\begin{equation}
\nonumber X^0=(C(w)\,z, b(w)),  \quad  X^1=(a(w)+A(w)(z,z), 2\,i\,< z, \bar{a}(w) >).
\end{equation}
Таким образом решения (\ref{88a}) -- это, в точности, четная компонента $X \in {\rm aut} \, \mathcal{Q}$ --
$X^0 \in {\rm aut}^0 \, \mathcal{Q}$, а решения (\ref{88b}) -- это нечетная компонента $X \in {\rm aut} \, \mathcal{Q}$ --
$X ^1\in {\rm aut}^1 \, \mathcal{Q}$. Отметим здесь, что  ${\rm aut}^0 \, \mathcal{Q}$, в отличие от ${\rm aut}^1 \, \mathcal{Q}$ -- это подалгебра ${\rm aut} \, \mathcal{Q}$.

\vspace{3ex}

 Как хорошо известно, подалгебра $g_{-}=g_{-2}+g_{-1}$ для любой невырожденной квадрики -- фундаментальна по Танаке (\cite{SS18}, \cite{G18}).  Поэтому если для некоторого $j$ сответствующая компонента обратилась в ноль $g_j=0$, то
$g_J=0$ для всех $J>j$.  В связи с этим дадим определение.

\vspace{3ex}
 {\bf Определение 2:} Назовем невырожденную $CR$-квадрику $\mathcal{Q}$ {\it исключительной}, если  для нее $g_3 \neq 0$.
\vspace{3ex}
Ясно, что квадрика не является исключительной тогда и только тогда, когда
$${\rm aut} \, \mathcal{Q}= g_{-2}+g_{-1}+g_{0}+g_{1}+g_{2}.$$
Примеры Мейлан и Грегоровича -- это примеры исключительных квадрик.

Условие отсутствия компонент веса больше, чем два в терминах $(a,A,b,C)$ -- это четыре условия
$$ {\rm deg}_u \, a \leq 1, \; {\rm deg}_u \, A =0, \;{\rm deg}_u \, C \leq 1, \;{\rm deg}_u \, b \leq 2.$$
Эти условия гарантируют, что четная компонента обрывается на втором весе, а нечетная - на первом.
Но поскольку критерием исключительнсти является условие на третью компоненту, то мы можем ограничится
рассмотрением только первой системы (\ref{88a}). Заметим также, что первое соотношение  (\ref{88a})
при фиксированном $a$ однозначно разрешимо относительно $A$ (условие невырожденности). Поэтому параметр $A$
можно исключить.

\vspace{3ex}

Через $z \cdot \bar{\zeta}$ будем обозначать стандартную эрмитову форму
$$(z \cdot \bar{\zeta})=z_1 \bar{\zeta}_1 + \dots +z_n  \bar{\zeta}_n .$$
Тогда произвольная эрмитова форма может быть записана как $ (H z \cdot \bar{z})$, где $H$ -- эрмитова матрица, а $z$ понимаем как столбец. Обозначим через $H_j z \cdot \bar{\zeta}$ -- $j$-ю координату вектор-формы $<z, \bar{z}>$.

\vspace{3ex}

Пусть $B(\bar{z})=(B_1 \cdot \bar{z},\dots,B_k  \cdot \bar{z})$ -- набор из $k$ линейных форм. Рассмотрим уравнение
\begin{equation}\label{B}
  <x,\bar{z}>=B(\bar{z}),     \quad x \in {\bf C}^n, \quad B(\bar{z})=(B_1 \cdot \bar{z},\dots,B_k  \cdot \bar{z}) \in {\bf C}^{kn},
\end{equation}
т.е. подразумевается, что при фиксированной правой части мы ищем вектор $x \in {\bf C}^n$, т.ч. уравнение превращается в тождество по $z$. Такое уравнение уравнение можно записать в виде $ H_j x =B_j, \; j=1,\dots,k$.
Рассмотрим далее отображение
 $$ \nu: {\bf C}^n  \rightarrow {\bf C}^{kn},   \quad \nu(z)= (H_1 z, \dots, H_k z).$$
Обозначим его образ через $\mathcal{L}'$, через $\mathcal{L}''$ обозначим прямое дополнение образа до всего пространства.
Пусть   $\pi$ -- проектор ${\bf C}^{kn}$ на $\mathcal{L}''$ вдоль $\mathcal{L}'$ .
Поскольку формы $(H_1, \dots,H_k)$ не имеют общего ядра, то $\nu$ -- это изоморофизм ${\bf C}^{n}$ и $\mathcal{L}'$, на $\mathcal{L}'$
определено обратное отображение $\nu^{-1}$. Рассуждения этого абзаца можно подытожить следующим образом.

\vspace{3ex}

{\bf Лемма 3 :} (a) Уравнение  (\ref{B})  разрешимо относительно $x$ тогда и только тогда, когда $\pi(B)=0$.\\
(b) Если $\pi(B)=0$, то единственное решение  (\ref{B}) имеет вид $x=\nu^{-1}(B)$.

\vspace{3ex}

Применяя лемму 3 к первому соотношению   (\ref{88a}), получаем

\vspace{3ex}

{\bf Лемма 4:}  (a) Уравнение
  $<A(u)(z,z),\bar{z} >=2 \,i \, < z, \Delta \bar{a}(u) >$
равносильно паре соотношений
\begin{equation}\label{Aa}
  \pi(< z, \Delta \bar{a}(u) > )=0,  \quad A(z,z)= 2 \, i \, \nu^{-1}(< z, \Delta \bar{a}(u) >).
\end{equation}
(b)  Если ${\rm deg}_u \, a = d$, то ${\rm deg}_u \, A \leq d-1$.

\vspace{3ex}

Прежде чем сформулировать полученную теорему, дадим определение.   Невырожденная квадрика называеется {\it жесткой},
если $g_{+}=0$. В силу фундаментальности это равносильно тому, что $g_1=0$.  А это, в свою очередь, равносильно тому, что
уравнение $\pi ( <z, \Delta \bar{a}(u)>)=0$ не имеет  ненулевых решений линейных по $u$.

Квадрику можно понимать как точку пространства наборов из $k$ эрмитовых форм от $n$-мерного переменного. Такие наборы -- это вещественное линейное пространство $\mathcal{H}$ размерности $k\,n^2$. Обозначим невырожденные квадрики через $\mathcal{H}'$.
Ясно, что вырожденные квадрики -- это алгебраическое подмножество $\mathcal{H}$ (задается полиномиальными условиями на координаты). Таким образом $\mathcal{H}'$, его дополнение, -- полуалгебраическое (равенства и неравенства на полиномы). Обозначим через $\mathcal{H}_l$ множество невырожденных квадрик,  для которых $g_{+}$ имеет длину не меньше $l$. Ясно, что не жесткие квадрики -- это
$\mathcal{H}_1$, а исключительные -- это $\mathcal{H}_3$. Из нашего описания компонент алгебры сразу следует, что $\mathcal{H}_l$ -- полуалгебраическое множество для всех $l$.
$$\mathcal{H}'  \supseteq \mathcal{H}_1 \supseteq  \mathcal{H}_2 \supseteq  \mathcal{H}_3 \supseteq \dots $$

Приведем два критерия:  исключительности и не жесткости.

\vspace{3ex}

{\bf Теорема 5:} Невырожденная квадрика $\mathcal{Q}$\\
(a) исключительна тогда и только тогда, когда существует ненулевая квадратичная вектор-форма
 $a(u,u)=(a_1(u,u),\dots,a_n(u,u)),$  удовлетворяющая двум соотношениям
 $$\pi(<z,\Delta \bar{a}(u,u)>)=0, \quad <z,\Delta^2 \bar{a}(u,u)>=0.$$
 (b) не жесткая тогда и только тогда, когда существует ненулевая линейная вектор-форма
 $a(u)=(a_1(u),\dots,a_n(u)),$  удовлетворяющая соотношению
 $$\pi(<z,\Delta \bar{a}(u)>)=0.$$

\vspace{3ex}

{\bf 2.  $CR$-квадрики типов $(3,3)$ и $(3,4)$}

\vspace{3ex}

Сформулируем два условия, которым может удовлетворять эрмитова вектор-форма $<z,\bar{z}>$ квадрики $\mathcal{Q}$.

(I)  Множество $\{z \in \mathbf{C}^n : {\rm rank}(H_1 z, \dots, H_k z)=n\}$ -- непусто (и, следовательно, открыто и плотно) .\\
(II) Образ отображения из $\mathbf{C}^{2n}$  в $\mathbf{C}^{k}$   $ \; (p,q) \rightarrow \; <p,q> \;$ содержит внутренние точки.

\vspace{3ex}

{\bf Теорема 6:} Если невырожденная квадрика удовлетворяет условиям (I) и (II), то она не является исключительной.\\
{\it Доказательство:}  Имеем $<z,\Delta^2 \bar{a}>=0$.  Из условия (I) следует, что $\Delta^2 a=d^2 a(<z,\bar{z}>,<z,\bar{z}>)=0$.
Окомплексим полученное равенство ( т.е. пусть $z$ и $\zeta=\bar{z}$ -- независимые переменные). Теперь из условия (II) следует, что в окрестности некоторого значения $z$ дифференциалы $du=<z,\bar{z}>$ могут рассматриваться как независимые. Это означает, что  $d^2 a=0$ и, следовательно,  ${\rm deg}_u \, a \leq 1$. Утверждение доказано.

Заметим, что доказательство утверждения основано только на втором соотношении критерия и не использует первого.

\vspace{3ex}

Зафиксируем  $ n > 0$. В силу условия линейной независимости координатных эрмитовых форм невырожденные квадрики коразмерности $k$ возможны лишь в диапазоне  $1 \leq k \leq n^2$.  При этом  в диапазоне $2 \leq k \leq n^2-2$ у квадрики общего положения $g_{+}=0$ (нет полей положительного веса). Критерием такой $"$жесткости$"$ является условие $g_1=0$.  В работе Н.Палинчак  \cite{NP} классифицированы с точностью до голоморфной эквивалентности все квадрики для $n=k=3$ с $g_1 \neq 0$. Это набор из восьми квадрик. Аналогичный список из девяти квадрик для $n=3 \; k=4$ был составлен С.Анисовой \cite{SA}. В силу фундаментальности из того, что $g_1=0$ следует обращение в ноль и последующих компонент. Поэтому исключительные квадрики данных типов, если они существуют, должны содержатся в списках Палинчак и Анисовой.

\vspace{3ex}

{\bf Теорема 7:} (a) Исключительных квадрик типа $(3,3)$  -- не существует.\\
(b) Исключительных квадрик типа $(n,4)$  при $n \leq 3$  -- не существует.\\
{\it Доказательство:} Непосредственно убеждаемся, что все квадрики из обоих списков удовлетворяют условию (II). Что касается условия (I),  то ему удовлетворяют все квадрики списка Анисовой и все квадрики списка Палинчак за одним исключением. Это квадрика, обозначенная там как $Q_5$. Определяющие  ее формы имеют вид.
$$2 \, {\rm Re} z_1\, \bar{z}_3,  \quad  2 \, {\rm Re} z_2\, \bar{z}_3,  \quad2 \, {\rm Im} z_1\, \bar{z}_3$$
Пространство  $B^{\bot}=\{z \in {\bf C}^3: <z,\bar{B}>=0\}$ задается двумя независимыми соотношениями
$z_2 \, \bar{B}_3+z_3 \, \bar{B}_2=z_1 \, \bar{B}_2-z_2 \, \bar{B}_1=0$. Поэтому решение имеет вид
$$B^{\bot}=\{z_1=\lambda \, \bar{B}_1, \;z_2=\lambda \, \bar{B}_2, \; z_3=-\lambda \, \bar{B}_3\}$$
Итак,  если $a(u)=(a_1(u),a_2(u),a_3(u))$ удовлетворяет соотношению $<\Delta^2  a, \bar{z}>=0$, то $\Delta^2 a_j(u)$ делится
на $\bar{z}_j$. Таким образом ограничение $\Delta^2  a_1$ на плоскость $\bar{z}_1=0$ равно нулю. При этом ограничения
координатных эрмитовых форм равны $(z_1 \, \bar{z}_3, \, z_2 \, \bar{z}_3+z_3 \, \bar{z}_2, -i \, z_1 \, \bar{z}_2)$.
Если рассмотреть соответствующее отображение $ {\bf C}^5$ на $ {\bf C}^3$, то мы видим, что оно имеет ранг три.
Это позволяет утверждать , что $d^2 a_1=0$ и степень $a_1$ не превосходит единицы. Аналогичные рассуждения дают ту же оценку степени для $a_2$ и  $a_3$. Таким образом для всех квадрик кроме $Q_5$ из списка Палинчак неисключительность  следует из утверждения 5, а для этой квадрики из нашего рассуждения. Это доказывает пункт (a) и отсутствие исключительных квадрик типа (3,4) . Квадрика типа (1,4) не может быть невырожденной.  Невырожденная квадрика типа (2,4) только одна. Это последняя квадрика, заданная базисом пространства эрмитовых форм на $ {\bf C}^2$.  Легко проверяется, что она не является исключительной (см.также утверждение 20) . Это доказывает (b). Теорема доказана.

\vspace{3ex}

{\bf 3. Квадрики малых коразмерностей}

\vspace{3ex}

Минимальная коразмерность, для которой известен пример исключительной квадрики это $k=4$.
С коразмерностями $k$ равными одному и двум дела обстоят следующим образом. То, что не существует исключительных квадрик в коразмерности один непосредственно следует из работы \cite{CM}, в \cite{BM20} было приведено доказательство для $k=2$.

В данной работе мы доказываем, что в коразмерности  $k=3$ исключительных квадрик также нет .
Таким образом,  коразмерность четыре -- это минимальная коразмерность, в которой возможно существование исключительных квадрик.

\vspace{3ex}

Пусть $f=\sum_{0}^{\infty} f_j,  \; g=\sum_{0}^{\infty} g_j$ -- разложения $f$ и $g$ в сумму весовых компонент. Тогда $X_m$, т.е. $m$-я весовая компонента поля (\ref{f}) имеет вид
\begin{eqnarray} \label{tc}
\nonumber X_m=2 \, {\rm Re} \, \left( f_{m+1}(z,w) \,\frac{\partial}{\partial z}
+ g_{m+2}(z,w) \,\frac{\partial}{\partial w}  \right),  \qquad\ \mbox{где} \qquad  \qquad\qquad\\
2 \, {\rm Re} \, \left( i \, g_{m+2} (z,u+i\,<z,\bar{z}>) + 2 \, <f_{m+1}(z,u+i\,<z,\bar{z}>) ,\bar{z}>  \right)=0.
\end{eqnarray}

\vspace{3ex}

В предположении невырожденности $\mathcal{Q}$  рассмотрим младшие компоненты алгебры.  При этом мы будем записывать возникающие здесь компоненты $f$ и $ g$ в виде полилинейных форм, причем будем предполагать, что они симметричны внутри каждой группы переменных $z$ или $w$. Имеем\\
{\it Вес (-2).}  $\quad f_{-1}=0, \; g_0=q,  \quad X_{-2}=(0,q)$.  Из условия (\ref{tc}) получаем, что $q \in \mathbf{R}^k$.\\
{\it Вес (-1).}  $\quad f_{0}=p, \; g_1=l(z)$,  где $p \in \mathbf{C}^n, \; l(z)$ -- линейная форма. Из условия (\ref{tc}) получаем, что
$l(z)=2\, i\, <z,\bar{p}>$. Таким образом  $X_{-1}=(p,\; 2\, i\, <z,\bar{p}>)$. \\
{\it Вес 0.}  $\quad f_{1}=C \,z , \; g_2=\alpha(z,z)+\rho \, w$. Из условия (\ref{tc}) получаем, что
$\alpha(z,z)=0, \;  {\rm Im} \,( \rho \,u)=0, \quad  2 \, {\rm Re}  <C \,z ,\bar{z}> =\rho \,<z ,\bar{z}> $.  Таким образом  $X_{0}=(Cz,\; \rho \,w)$ с найденными условиями. \\
{\it Вес 1.}  $\quad f_{2}=a \,w +A(z,z), \; g_3=\alpha(z,z,z)+\beta(z) \, w$. Из условия (\ref{tc}) получаем, что
$\alpha(z,z,z)=0, \; <A(z,z),\bar{z}>=2 \, i \, <z,\bar{a} <z,\bar{z}>>,  \; \beta(z) \, u = 2 \, i \, <z,\bar{a} \, u> $.
Получаем $X_1=(aw+A(z,z), \, 2 \, i \, <z,\bar{a} \, w>)$.\\
{\it Вес 2.}  $\quad f_{3}=B(w)z +b(z,z,z), \; g_4=\alpha(z,z,z,z)+\beta(z,z) \, w+r(w,w)$. Из условия (\ref{tc}) получаем, что
$\alpha(z,z,z,z)=0, \; \beta(z,z) \, u=0, \; b(z,z,z)=0$ ,  а также
${\rm Re} <B(u)z,\bar{z}>=r(<z,\bar{z}>,u) \; {\rm Im} <B(<z,\bar{z}>)z,\bar{z}>=0$.
Получаем $X_2=(B(w)\,z, \, r(w,w))$ с найденными условиями.\\
{\it Вес 3.}  $\quad f_{4}=d(w,w)+D(w)(z,z) +e(z,z,z,z), \; g_5=\alpha(z,z,z,z,z)+\beta(z,z,z) \, w+\gamma(z)(w,w)$. Из условия (\ref{tc}) получаем, что
$\alpha(z,z,z,z,z)=0,  \; \beta(z,z,z)(u)=0, \; e(z,z,z,z)=0$ ,  а также
\begin{eqnarray*}
\gamma(z)(u,u)=2\,i\,<z,\bar{d}(u,u)>, \\
<D(u)(z,z),\bar{z}>=4\,i\,<z,\bar{d}(<z,\bar{z}>,u)>, \quad  <D(<z,\bar{z}>)(z,z),\bar{z}>=0.
\end{eqnarray*}
Получаем $X_3=(d(w,w)+D(w)(z,z), \, 2\,i\,<z,\bar{d}(w,w)>)$ с найденными условиями.

Итак, получаем

\vspace{3ex}

{\bf Утверждение 8:}  (a) Если $\mathcal{Q}$ -- невырожденная квадрика, то младшие весовые компоненты ${\rm aut} \, \mathcal{Q}$ имеют вид
\begin{eqnarray*}
X_{-2}=(0,q),   \quad  q \in \mathbf{R}^k,
\qquad\qquad\qquad\qquad\qquad\qquad\qquad\qquad\qquad\qquad\qquad\qquad\qquad\qquad
\qquad\qquad\qquad\qquad\qquad\qquad\qquad\qquad\qquad\\
X_{-1}=(p,\; 2\, i\, <z,\bar{p}>), \quad p \in \mathbf{C}^n,\qquad\qquad\qquad\qquad\qquad\qquad\qquad\qquad\qquad\qquad\qquad\qquad\qquad\qquad\qquad\qquad
\qquad\qquad\qquad\qquad\qquad\qquad\qquad\\
X_{0}=(Cz,\; \rho \,w),  \; C \in gl(n,\mathbf{C}), \qquad  \rho \in gl(k,\mathbf{R}),  \quad
  2 \, {\rm Re}  <C \,z ,\bar{z}> =\rho \,<z ,\bar{z}>,\qquad\qquad\qquad\qquad\qquad\qquad\qquad\qquad
  \qquad\qquad\qquad\qquad\qquad\qquad\qquad\qquad\qquad\qquad\\
X_1=(aw+A(z,z), \, 2 \, i \, <z,\bar{a} \, w>), \qquad
 <A(z,z),\bar{z}>=2 \, i \, <z,\bar{a} <z,\bar{z}>>,\qquad\qquad\qquad\qquad\qquad\qquad\qquad\qquad
 \qquad\qquad\qquad\qquad\qquad\qquad\qquad\qquad\qquad\qquad\qquad\\
X_2=(B(w)\,z, \, r(w,w)), \quad
{\rm Re} <B(u)z,\bar{z}>=r(<z,\bar{z}>,u), \;\; {\rm Im} <B(<z,\bar{z}>)z,\bar{z}>=0, \qquad\qquad\qquad\qquad\qquad\qquad\qquad\qquad
\qquad\qquad\qquad\qquad\qquad\qquad\qquad\qquad\qquad
\end{eqnarray*}
\begin{eqnarray} \label{g3}
\nonumber
X_3=(d(w,w)+D(w)(z,z), \, 2\,i\,<z,\bar{d}(w,w)>), \qquad\qquad\qquad\qquad\qquad\qquad\\
 <D(u)(z,z),\bar{z}>=4\,i\,<z,\bar{d}(<z,\bar{z}>,u)>, \;  <D(<z,\bar{z}>)(z,z),\bar{z}>=0 \qquad\qquad
 \end{eqnarray}
(b)  $\mathcal{Q}$ - не является исключительной тогда и только тогда, когда единственным решением  (\ref{g3})
является $D(u)=0$   и $d(u,u)=0$.
\vspace{3ex}

{\bf Лемма 9 :} (a) Если $\mathcal{Q}$ -- невырожденная квадрика, то
$$D(<z,\bar{p}>)(z,z)=0.$$
(b)  $\mathcal{Q}$ - не является исключительной тогда и только тогда, когда единственным решением
\begin{equation}\label{gg3}
 D(<z,\bar{p}>)(z,z),\bar{z}=0,  \quad <D(u)(z,z),\bar{z}>=4\,i\,<z,\bar{d}(<z,\bar{z}>,u)>.
\end{equation}
является $D(u)=0$   и $d(u,u)=0$.\\
{\it Доказательство:} Докажем (a). Непосредственное вычисление коммутатора дает
$[X_3,X_{-1}]=(F(z,w),G(z,w))$, где
\begin{eqnarray*}
F=-2\, d(2i<z,\bar{p}>,w)-D(2i<z,\bar{p}>)(z,z)-2D(w)(z,p),\qquad\qquad\qquad\qquad\\
G=(2\,i\,<d(w,w)+D(w)(z,z),\bar{p}>-2\, i\, <p,\bar{d}(w,w)>-4\, i\, <z,\bar{d}(2i<z,\bar{p}>,w)>.
\end{eqnarray*}
Коммутатор поля веса 3 с полем веса (-1) это элемент $g_2$. Нам известно, что $z$-координата такого поля не содержит кубических форм от $z$, поэтому  $D(<z,\bar{p}>)(z,z)=0$. Отсюда сразу следует (b). Лемма доказана.

\vspace{3ex}

Пусть $D(u)(z,z)=D_1(z,z) \, u_1+ \dots + D_k(z,z) \, u_k$, где $D_j(z,z)=(D_j^1(z,z), \dots, D_j^n(z,z))$ - векторзначная квадратичная форма. Тогда первое соотношение  (\ref{gg3}) принимает вид
\begin{equation}\label{D1}
(H_1 z \cdot \bar{p}) \, D_1(z,z)+ \dots + (H_k z \cdot \bar{p}) \, D_k(z,z)=0
\end{equation}
При этом  $\nu$-я координата (\ref{D1}) выглядит так
\begin{equation}\label{D2}
(H_1 z \cdot \bar{p}) \, D^{\nu}_1(z,z)+ \dots + (H_k z \cdot \bar{p}) \, D^{\nu}_k(z,z)=0
\end{equation}
Убирая свертку с независимым параметром $\bar{p}$, получаем векторное соотношение вида
 \begin{equation}\label{D3}
  D^{\nu}_1(z,z) \, H_1 z + \dots +  D^{\nu}_k(z,z) \, H_k z=0
 \end{equation}

 \vspace{3ex}

Квадратичной вектор-форме $d(u,u)$ однозначно соответствует симметрическая билинейная вектор-форма
\begin{eqnarray}
\nonumber
d(u,U)=\sum_{\alpha\beta} \, d_{\alpha\beta}\, u_{\alpha} \, U_{\beta} = \sum_{\beta} \, \delta_{\beta}(u) \, U_{\beta}, \quad
\delta_{\beta}(u)=\sum_{\alpha} \, d_{\alpha\beta}\, u_{\alpha},
 \quad d_{\alpha\beta} \in \mathbf{C}^n.
\end{eqnarray}

Если $D(u)=0$, то  условие на форму $d(u,u)$ можно записать так
$$
 <  z,\bar{\delta}_j(<z,\bar{z}>)>=0 ,   \quad j=1,\dots,k.
$$
Поскольку все коэффициенты формы $d$ удовлетворяют одному и тому же соотношению запишем его в свободном от указания на индекс виде.
\begin{equation}\label{d}
 <  z,\bar{\delta}(<z,\bar{z}>)>=0, \quad   \delta(u)=u_1 \, \alpha_1+\dots+u_k \, \alpha_k,
\end{equation}
Таким образом, чтобы описать $g_3$ при условии $D(z,z)=0$ нужно описать все наборы $(\delta_1(u),\dots,\delta_k(u))$  линейных $\mathbf{C}^n$-значных форм , таких, что каждая удовлетворяет соотношению (\ref{d}), а вся совокупность удовлетворет условию симметричности билинейной формы $d(u,U)$, а именно $d_{\alpha\beta}=d_{\beta\alpha}$.

\vspace{3ex}

Прежде, чем перейти к рассмотрению интересующего нас случая коразмерности три,  рассмотрим в нашем контексте случаи $k=1$ и $k=2$.\\

{\it Случай k=1.}
 Если $k=1$, то первое соотношение (\ref{gg3}) принимает вид $D(z,z)<z,\bar{p}>=0$, откуда сразу получаем, что $D=0$ После
 второе  соотношение принимает вид
$<d(<z,\bar{z}>,u), \bar{z}>=<\alpha, \bar{z}><z,\bar{z}>u=0$. Откуда сразу получаем $d=0$ и $g_3=0$.

 \vspace{3ex}

{\it Случай k=2.}  Положим
$$Pz =H_1 z, \; Qz=H_2z, \quad  Pz=(p_1(z),...,p_n(z)), \; Qz=(q_1(z),...,q_n(z)).$$
Тогда (\ref{D3}) принимает вид
  $$ l(z,z) \, P z +  m(z,z) \, Q z=0,    \; l(z,z)=D^{\nu}_1(z,z), \; m(z,z)=D^{\nu}_2(z,z)$$
Пусть найдется такое $\mu$, что $p_{\mu}(z)=\alpha(z)$ не пропорционально   $q_{\mu}=\beta(z)$. Тогда получаем, что
 $l(z,z)=\lambda(z) \, \beta(z), \; m(z,z)=-\lambda(z) \, \alpha(z)$. Если $\lambda(z) \neq0$, то получаем
 $\beta(z) \,H_1 z=\alpha(z) \, H_2 \,z$, откуда следует, что $H_1 z =\alpha(z) \, A, \, H_2 z =\beta(z) \, A$.

 Сформулируем лемму, которая понадобится нам здесь и в дальнейшем.

 \vspace{3ex}

{\bf Лемма 10:} (a) Если эрмитова матрица $H$ имеет ранг 1, т.е. $H z =l(z) A, \; l \neq 0,  A \neq 0$, то $l(z)=\rho \,(z \cdot \bar{A}), \, \rho \neq 0$, а соответствующая эрмитова форма имеет вид
$$(Hz \cdot \bar{z})=\rho \, | (z \cdot \bar{A})|^2, \; \rho \in \mathbf{R}.$$
(b) (a) Если эрмитова матрица $H$ имеет ранг 2, т.е. $H z =l(z) A+m(z) B,$  причем $l$ и $m$ -- линейно независимы ${\rm rank}( A,B)=2$, то $l(z)=\rho \,(z \cdot \bar{A}),  \; m(z)=\tau \,(z \cdot \bar{B}),$ параметры $\rho$ и $\tau$ -- вещественны и не
равны нулю, а соответствующая эрмитова форма имеет вид
$$(Hz \cdot \bar{z})=\rho \, | (z \cdot \bar{A})|^2+ \tau \, | (z \cdot \bar{B})|^2.$$
{\it Доказательство:} Ранг эрмитовой формы инвариантен относительно невырожденных линейных замен. Приведем формы к диагональному виду.  Число ненулевых коэффициентов в первом случае равно одному, во втором -- двум. Возвращаясь к исходным координатам получаем оба утверждения леммы.

  \vspace{3ex}
В соответствии с леммой
$$(P z \cdot \bar{z}) = \tau_1 |(z \cdot \bar{A})|^2, \quad (Q z \cdot \bar{z}) = \tau_2 |(z \cdot \bar{A})|^2,$$
что противоречит условию невырожденности, поэтому $D=0$.
  \vspace{2ex}

Запишем (\ref{d}), получим
\begin{eqnarray} \label{d2}
(P z \cdot \bar{z})\,<z,\bar{\alpha}_1> +(Q z \cdot \bar{z})\,<z,\bar{\alpha}_2>=0
\end{eqnarray}
Пара векторов $(\alpha_1, \alpha_2)$ может иметь ранг ноль, один или два. Ноль -- это  значит, что $\delta=0$. Пусть -- один, тогда вектора пропорциональны, пусть $\alpha_2=\lambda \, \alpha_1,  \; \alpha_1 \neq0$. Из (\ref{d2}) получаем
\begin{eqnarray*}
((P z \cdot \bar{z}) +\bar{\lambda} \,(Q z \cdot \bar{z}))\,<z,\bar{\alpha}_1>=0
\end{eqnarray*}
Откуда следует линейная зависимость $P$ и $Q$,  что противоречит условию невырожденности пары $(P,Q)$.

Пусть теперь $(\alpha_1, \alpha_2)$ -- линейно независимы ($n \geq 2$).  Выберем базис в $\mathbf{C}^n$ вида
$(\alpha_1, \alpha_2, \dots)$. Пусть $p_j(z)$ -- это $j$-я координата $Pz$ в этом базисе, а  $q_j(z)$ -- это $j$-я координата $Qz$.
Тогда  (\ref{d2}) примет  вид
\begin{eqnarray}\label{d22}
 p_1(z)\, P z +p_2(z)\,Qz=0, \quad  q_1(z)\, P z +q_2(z)\,Qz=0
\end{eqnarray}
Выпишем первую координату первого равенства и вторую -- второго. Имеем
$$(p_1(z))^2+p_2(z) \,q_1(z)=0,  \quad (q_2(z))^2+p_2(z) \, q_1(z)=0$$
Если $p_1=0$, то, как следует из первого равенства  $p_2=0$ (иначе $q_1=0$ и $\alpha_1$ попадает в общее ядро).
Теперь из второго равенства следует, что $q_2=0$ и $\alpha_2$ попадает в общее ядро. Противоречие.
Значит $p_1$ и $q_2$ не равны нулю, но тогда из  (\ref{d22}) и невырожденности следует, что   $p_2$ и $q_1$ также не равны нулю.
Тогда, применяя лемму 3 к первому равенству  (\ref{d22}) , видим, что возможен только вариант (b.2), т.е.
$Pz=p_2(z)\,A, \; Qz= -p_1(z)\,A$. В соответствии с леммой 6 формы  $(Pz \cdot \bar{z})$ и $(Qz \cdot \bar{z})$ пропорциональны  $| (z \cdot \bar{A})|^2$, что противоречит условию невырожденности.

Таким образом для $k=2$  компонента $g_3=0$.

\vspace{7ex}

 В связи с изучением (\ref{gg3}) для малых $k$ cформулируем еще несколько очевидных вспомогательных утверждений.

\vspace{3ex}

{\bf Лемма 11:}  Пусть $m(z)\,P\,z=l(z)\,Q\,z$,  где $P$ и $Q$ -- квадратные матрицы, а $m$ и $l$ -- скалярные линейные формы,   тогда реализуется  одна из следующих возможностей:\\
(a)  $m(z)\,P\,z=l(z)\,Q\,z=0$ --  в каждой паре хотя бы один из сомножителей равен нулю (4 варианта).\\
(b.1) Линейная зависимость.  $P \, z= \lambda \, Q \, z \neq 0, \;\; m(z)=\lambda \, l(z) \neq 0$.\\
(b.2)  Линейная независимость. $P\,z =l(z) \, A, \;\; Q\,z =m(z) \, A, \;\; A \neq 0$ , формы $l$ и $m$ не пропорциональны.\\

\vspace{3ex}
{\bf Лемма 12:}
Пусть $m(z)\,A=l(z)\,B$,  где $A$ и $B$ -- вектора, а $m$ и $l$ -- скалярные линейные формы,   тогда реализуется  одна из следующих возможностей:\\
(a)  $m(z)\,A=l(z)\,B=0$ --  в каждой паре хотя бы один из сомножителей равен нулю (4 варианта).\\
(b) Линейная зависимость.  $B= \lambda \, A \neq 0, \;\; m(z)=\lambda \, l(z) \neq 0$.\\

\vspace{3ex}
Пусть $I(l_1(z),l_2(z))$ -- идеал в кольце полиномов от $z$, порожденный двумя линейными формами $l_1$ и $l_2$.

{\bf Лемма 13:}  (a)  $l(z) \in I(l_1(z),l_2(z))$ тогда и только тогда, когда $l(z)=\lambda_1 \, l_1(z)+\lambda_2 \, l_2(z)$.\\
(b) Идеал $I(l_1(z),l_2(z))$ -- прост, т.е. если $p(z) \, q(z) \in I(l_1(z),l_2(z))$, то либо $p(z) \in  I(l_1(z),l_2(z))$, либо $q(z) \in  I(l_1(z),l_2(z))$.

\vspace{3ex}

{\bf Перейдем к рассмотрению коразмерности } $\bf k=3$. \\

\vspace{3ex}

Пусть
$$<z,\bar{z}>=((P z \cdot \bar{z}), \,( Q z \cdot \bar{z}), ( R z \cdot \bar{z}))$$
Соотношение (\ref{D3}) примет вид
\begin{equation}\label{DD3}
l(z,z) \, P z+m(z,z) \, Q z+n(z,z) \, R z=0
\end{equation}

\vspace{3ex}

{\bf Лемма 14:} Соотношение (\ref{DD3}) невозможно для линейно зависимых $(l,m,n)$.\\
{\it Доказательство:} Пусть ${\rm rank} (l,m,n)=1$, т.е.
$$l(z,z)=\lambda \, \varphi(z,z), \; m(z,z)=\mu \, \varphi(z,z), \; n(z,z)=\nu \, \varphi(z,z),   \varphi(z,z) \neq 0.$$
Тогда $\lambda \, P z+\mu \, Q z+\nu \, R z=0$. Противоречие.\\
Пусть ${\rm rank} (l,m,n)=2$, т.е.  (для определенности)  $n(z,z) = \lambda \, l(z,z)+\mu \, m(z,z)$, где $l$ и $m$ -- непропорциональны.
Имеем
$$l(z,z) (P z + \lambda \, R z)+m(z,z) (Qz+\mu \,  R z)=0.$$
Пусть $l$ и $m$ -- не взаимно просты, т.е.
$$l(z,z)=\alpha(z) \, \gamma(z), \; m(z,z)=\beta(z) \, \gamma(z),$$
где $\alpha$ и $\beta$ -- непропорциональны $\gamma \neq 0$. Получаем
$$ \alpha(z) (P z + \lambda \, R z)+\beta(z) (Qz+\mu \,  R z)=0.$$
Откуда следует (лемма 11), что
$$(P z + \lambda \, R z)=\beta(z) \, A, \; (Qz+\mu \,  R z)=-\alpha(z) \, A, \;    A\in \mathbf{C}^n.$$
Если $\lambda$ и $\mu$ -- вещественны, то матрицы, которые стоят в левых частях этих равенств -- эрмитовы и по лемме 10
$\alpha(z)$ и $\beta(z)$ пропорциональны $(z \cdot \bar{A})$. Противоречие.\\
Если хотя бы одно из чисел $\lambda$ и $\mu$ -- не вещественно, например $\lambda$, то, отделяя в
$(P z \cdot \bar{z} + \lambda \, R z\cdot \bar{z})=\beta(z) \, A\cdot \bar{z} $ мнимую часть получаем
$$({\rm Im} \lambda ) \,  R z\cdot \bar{z} ={\rm Im}(\beta(z) \, A\cdot \bar{z})$$
Откуда по лемме 10 получаем, что $\beta(z)$ пропорциональна $z\cdot \bar{A}$,  откуда следует, что все три эрмитовы формы
$(P z \cdot \bar{z}, \; Q z \cdot \bar{z}, \;  R z\cdot \bar{z})$ пропорциональны $|z\cdot \bar{A}|^2$. Противоречие.\\
Пусть $l$ и $m$ -- взаимно просты,  тогда $(P z + \lambda \, R z)=(Qz+\mu \,  R z)=0$. Противоречие. Лемма доказана.

\vspace{3ex}

Далее полагаем, что $(l(z,z),m(z,z),n(z,z))$ -- линейно независимы.  Рассмотрим идеалы, порожденные парами квадратичных форм
$$I_1=(m,n), \; I_2=(l,n), \; I_3=(l,m)$$
и множества их нулей $V_1,V_2,V_3$.  Линейная независимость форм означает, что
\begin{equation}\label{nn}
 l(z,z) \notin I_1, \quad m(z,z) \notin I_2, \quad n(z,z) \notin I_3.
\end{equation}
Рассмотрим множества нулей этих идеалов
$$V_1=\{m(z,z)=n(z,z)=0\}, \;V_2=\{l(z,z)=n(z,z)=0\}, \;V_3=\{l(z,z)=m(z,z)=0\}$$
\vspace{3ex}

{\bf Лемма 15:}  Пусть размерность пространства $ {\bf C}^n$ не меньше 4 и $l$ не есть тождественный ноль на $V_1$,
$m$ не есть тождественный ноль на $V_3$, $n$ не есть тождественный ноль на $V_3$. Тогда соотношение (\ref{DD3}) невозможно.
{\it Доказательство:} Тогда из  соотношения (\ref{DD3}) следует, что $Pz=0$  на непустом множестве
$$V'_1=\{z \in {\bf C}^n : m(z,z)=n(z,z)=0, \; l(z,z) \neq 0\}.$$
Поэтому  $Pz=0$ и на множестве $\{z \in {\bf C}^n : m(z,z)=n(z,z)=0\},$ которое лежит в замыкании $V_1$. В силу аналогичных рассуждений для $Qz$ и $Rz$ все три линейных оператора обращаются в ноль на
 $$\{z \in {\bf C}^n : m(z,z)=n(z,z)=l(z,z)= 0\},$$
 которое  имеет размерность не ниже $(n-3)$ и, поэтому содержит $z \neq 0$. Противоречие. Лемма доказана.

 \vspace{3ex}

Рассмотрим ситуацию, когда хотя бы одно из трех условий предыдущей леммы нарушено. Пусть, например, $l=0$ на $V_1$.
Учитывая, что $(l,m,n)$ -- линейно независимы  и $l \notin I_1$ из этого сразу следует, что $I_1$ не является радикальным.
Это означает, что среди линейных комбинаций $m$ и $n$ есть полный квадрат. Вместе  с условием обращения в ноль $l$
это означает, что от линейно независимых форм $(l,m,n)$ мы можем с помощью невырожденного линейного преобразования
перейти к их линейным комбинациям вида
$$ l'(z,z)=\lambda(z) \, \mu(z),  \; m'(z,z)=(\mu(z))^2, \;  n'(z,z)=n(z,z)$$
где линейные формы $\lambda$ и $\mu$ - непропорциональны.  Соотношение (\ref{DD3}) примет вид
\begin{equation}\label{DDD3}
l'(z,z) \, P' z+m'(z,z) \, Q' z+n'(z,z) \, R' z=\lambda(z) \, \mu(z) \, Pz+(\mu(z))^2 \, Qz+n(z,z) \,Rz=0
\end{equation}
где $(P'z,Q'z,R'z)$ получены из $(Pz,Qz,Rz)$ невырожденным линейным преобразованием.

 \vspace{3ex}

{\bf Лемма 16:}  Если  размерность пространства $ {\bf C}^n$ не меньше 4 и одно из трех условий леммы 15 нарушено, то выполнение  (\ref{DD3}) -- невозможно.
{\it Доказательство:} Пусть $n(z,z)$ делится на $\mu(z)$, т.е. $n(z,z)=\mu(z) \, \nu(z)$, причем $(\lambda(z),\mu(z),\nu(z))$ -- линейно независимы. Тогда (\ref{DDD3}) принимает вид
$$ \lambda(z) \, Pz+\mu(z) \, Qz+\nu(z) \,Rz=0.$$
Тогда, также как при доказательстве леммы 15 покажем, что все три оператора обращаются в ноль на подпространстве
$$ \lambda(z)= \mu(z) = \nu(z)=0$$
положительной размерности. Это означает наличие общего ядра у исходной тройки операторов. Противоречие.\\

    Пусть теперь $n(z,z)$ не делится на $\mu(z)$,  имеем
$$\mu(z) ( \lambda(z) \, Pz+\mu(z) \, Qz) + n(z,z) \, Rz=0.$$
 Откуда следует, что
 $$Rz=\mu(z) \, A , \quad   ( \lambda(z) \, Pz+\mu(z) \, Qz) + n(z,z) \, A=0, \quad A \neq 0.$$
Откуда получаем, что $n(z,z)=\alpha(z) \, \lambda(z) +\beta(z) \,\mu(z)$. Тогда  имеем
$$\lambda(z)(Pz+\alpha(z)\,A)+\mu(z)(Qz+\beta(z)\,A)=0.$$
Откуда, в свою очередь, получаем
$$(Pz+\alpha(z)\,A)=\mu(z) \,B,  \quad (Qz+\beta(z)\,A)=-\lambda(z) \,B.$$
Итого
$$Pz=-\alpha(z)\,A+\mu(z) \,B,  \quad Qz=-\beta(z)\,A)-\lambda(z) \,B, Rz=\mu(z) \, A.$$
Т.е. $(Pz,Qz,Rz)$ лежат в линейной оболочке двух постоянных векторов $(A,B)$.  От $(Pz,Qz,Rz)$
с помощью невырожденного линейного преобразования можно перейти к исходным эрмитовым матрицам,
которые фигурируют в соотношении   (\ref{DD3}). При этом векторам $(A,B)$ будут соответствовать
два простоянных вектора $(A',B')$.  В таком случае, как следует из леммы 10, все три эрмитовых формы
есть линейные комбинации $|(z \cdot \bar{A}'|^2$ и  $|(z \cdot \bar{B}'|^2$ и, следовательно, линейно зависимы.
Проттворечие.

Применяя лемму 14, лемму 15, лемму 16 и теорему 7 получаем следующее утверждение

{\bf Утверждение 17:}  (a) Пусть $k=3$.  Если $D(u)(z,z)$ удовлетворяет соотношению (\ref{D1}), то $D(u)(z,z)=0$.\\
(b) При этом $g_3$ состоит из полей вида $X=(d(w,w), \, 2\,i\,<z,\bar{d}(w,w)>)$, где  $<z,\bar{d}(<z,\bar{z}>,u)>=0$.

 \vspace{5ex}

Для $k=3$ форма $\delta$ имеет вид.

$$\delta(u)=\alpha_1 \, u_1+\alpha_2 \, u_2+\alpha_3 \, u_3.$$ Запишем (\ref{d}),
получим
\begin{eqnarray} \label{d2}
(P z \cdot \bar{z})\,<z,\bar{\alpha}_1> +(Q z \cdot \bar{z})\,<z,\bar{\alpha}_2>+(R z \cdot \bar{z})\,<z,\bar{\alpha}_3>=0
\end{eqnarray}
Тройка векторов $(\alpha_1, \alpha_2, \alpha_3)$ может иметь ранг ноль, один, два или три. Ноль -- это  значит, что $\delta=0$.

{\it Случай 1 (ранг один).}
Пусть ранг равен одному, т.е. $\alpha_{j}= \lambda_j \, \alpha, \; \alpha \in \mathbf{C}^n \setminus \{0\}, \; (\lambda_1,\lambda_2,\lambda_3) \neq 0$.
 Из (\ref{d2}) получаем
\begin{eqnarray*}
(\bar{\lambda}_1 \,(P z \cdot \bar{z}) +\bar{\lambda}_2 \,(Q z \cdot \bar{z})+\bar{\lambda}_3 \,(R z \cdot \bar{z}))\,<z,\bar{\alpha}>=0
\end{eqnarray*}
Откуда следует линейная зависимость $(P,Q,R)$,  что противоречит условию невырожденности этого набора.

{\it Случай 2 ( ранг два).}  Пусть ранг равен двум. Можем считать, что $\alpha_{3}= \bar{\lambda }\, \alpha_1+\bar{\mu} \, \alpha_2$, где вектора
$\alpha_1, \alpha_2$ - линейно независимы.
 Из (\ref{d2}) получаем
\begin{eqnarray*}
((P z \cdot \bar{z}) +\lambda \,(R z \cdot \bar{z})) \, <z,\bar{\alpha}_1> +((Q z \cdot \bar{z})+\mu \,(R z \cdot \bar{z}))\, <z,\bar{\alpha}_2>=0.
\end{eqnarray*}
Выберем в $\mathbf{C}^n$  базис вида $(\alpha_1, \alpha_2, \dots)$.
Записывая (\ref{d2}), получаем
\begin{eqnarray} \label{d32}
\nonumber
 p_1(z)\,((P z \cdot \bar{z}) +\lambda \,(R z \cdot \bar{z})) + p_2(z) \,((Q z \cdot \bar{z}) +\mu \,(R z \cdot \bar{z}))=0,\\
\nonumber
 q_1(z)\,((P z \cdot \bar{z}) +\lambda \,(R z \cdot \bar{z})) + q_2(z) \,((Q z \cdot \bar{z}) +\mu \,(R z \cdot \bar{z}))=0,\\
  r_1(z)\,((P z \cdot \bar{z}) +\lambda \,(R z \cdot \bar{z})) + r_2(z) \,((Q z \cdot \bar{z}) +\mu \,(R z \cdot \bar{z}))=0.
\end{eqnarray}
Из условия невырожденности набора $(P,Q,R)$ следует, что в каждой паре
$(p_1,p_2), (q_1,q_2),(r_1,r_2)$ либо обе формы равны нулю, либо -- непропорциональны. Пусть $r_1=r_2=0$.
Тогда записывая первую координату первого соотношения и вторую -- второго и применяя лемму 5, получаем, что
$p_1$ пропорциональна $p_2$,  а $q_1$ пропорциональна $q_2$. Таким образом все эти формы равны нулю, а это означает
что оба вектора $\alpha_1$ и $\alpha_2$ - содержатся в ядре $<z ,\bar{z}>$. Это противоречие означает, что $(r_1,r_2)$ -- непропорциональны.
Тогда из леммы 3 и условия невырожденности из третьего соотношения (\ref{d32}) получаем
\begin{equation}\label{d322}
P\,z +\lambda\,R\,z=r_2(z)\,A, \quad Q\,z +\mu\,R\,z=-r_1(z)\,A, \quad A\neq 0
\end{equation}
Подставляя (\ref{d322}) в первое и второе соотношения (\ref{d32}) и учитывая, что $A\neq0$, получаем
\begin{eqnarray*}
  (p_1(z)\,r_2(z)-p_2(z)\,r_1(z))\,(A \cdot \bar{z}))=0, \quad (q_1(z)\,r_2(z)-q_2(z)\,r_1(z))\,(A \cdot \bar{z}))=0.
\end{eqnarray*}
Учитывая непропорциональность $r_1$ и $r_2$, получаем
$$p_1(z)=\nu \, r_1(z), \; p_2(z)=\nu \, r_2(z), \; q_1(z)=\kappa \, r_1(z), \; q_2(z)=\kappa \, r_2(z).$$
Выписывая первую и вторую координаты (\ref{d322}), получаем
\begin{eqnarray*}
\nu \, r_1(z)+\lambda\,r_1(z)=a_1 \, r_2(z),   \quad  \nu \, r_2(z)+\lambda\,r_2(z)=a_2 \, r_2(z),\\
\kappa \, r_1(z)+\mu \,r_1(z)=a_1 \, r_1(z),   \quad  \kappa \, r_2(z)+\mu \, r_2(z)=a_2 \, r_1(z).
\end{eqnarray*}
Откуда следует, что
$$\nu=-\lambda, \; \kappa=-\mu, \; a_1=a_2=0.$$
Перепишем (\ref{d322}) в виде
\begin{equation}\label{d4}
(P z \cdot \bar{z}) +\lambda \,(R z \cdot \bar{z})=r_2(z)\,(A \cdot \bar{z}), \; (Q z \cdot \bar{z}) +\mu \,(R z \cdot \bar{z})=-r_1(z)\,(A \cdot \bar{z})
\end{equation}
Если коэффициенты $\lambda$ и $\mu$ -- вещественны, то из леммы 6 следует, что обе формы $r_1$ и $r_2$ пропорциональны
$(z \cdot \bar{A})$. Противоречие. Пусть, далее, лишь одно из них вещественно, скажем, $\mu$, а
$L= {\rm Im}\,\lambda \neq 0$.  Тогда из второго соотношения (\ref{d4}) получаем, что $r_1(z)=\tau\,(z \cdot \bar{A})$.
Отделяя мнимую часть первого соотношения (\ref{d4}) получаем, что
$$(R z \cdot \bar{z})=\frac{1}{L}\, {\rm Im}\,(r_2(z)\,(A \cdot \bar{z})).$$
Это равносильно тому,что
$$r_j(z)=\frac{1}{2\,i\,L}(r_2(z)\,a_j -( z\cdot \bar{A}) \, \bar{r}_2^j).$$
Для $j=2$ получаем
$$r_2(z)=\frac{\bar{r}_2^2}{2\,i\,L}\, ( z\cdot \bar{A}).$$
Откуда следует, что $r_1$ и $r_2$ пропорциональны. Противоречие.

Пусть теперь $L= {\rm Im}\,\lambda \neq 0$ и $M= {\rm Im}\,\mu \neq 0$.
Отделяя мнимую часть первого и второго соотношения (\ref{d4}) получаем, что
$$(R z \cdot \bar{z})=\frac{1}{L}\, {\rm Im}\,(r_2(z)\,(A \cdot \bar{z}))=-\frac{1}{M}\, {\rm Im}\,(r_1(z)\,(A \cdot \bar{z})).$$
Откуда следует, что
$$(z \cdot \bar{A})\,(L\,r_1(z)+M \,r_2(z))=(z \cdot \bar{A})\,(L\,\overline{r_1(z)}+M \,\overline{r_2(z)}).$$
Откуда следует, что
$$r_2(z)=-\frac{L}{M} \,r_1(z)  + \nu \, ( z\cdot \bar{A}). $$
Отделяя теперь вещественные части соотношений (\ref{d4}), подставляя туда полученные значения для $(R z \cdot \bar{z})$ и $r_2(z)$,
мы получаем выражения для $(P z \cdot \bar{z})$ и $(Q z \cdot \bar{z})$. В результате мы видим, что все три эрмитовых формы пропорциональны ${\rm Im}\,(r_1(z)\,(A \cdot \bar{z}))$. Противоречие.

{\it Случай 3 (ранг три).}  Будем предполагать, что $n \geq 3$.  То, что не существует контрпримера $CR$-размерности $n \leq 2$ -- хорошо известно. Выберем в $\mathbf{C}^n$  базис вида $(\alpha_1, \alpha_2,  \alpha_3,\dots)$.  Соотношение (\ref{d2}) примет вид
\begin{eqnarray} \label{d3}
\nonumber
 p_1(z)\,(P z \cdot \bar{z}) + p_2(z) \,(Q z \cdot \bar{z}) +p_3(z) \,(R z \cdot \bar{z})=0,\\
\nonumber
q_1(z)\,(P z \cdot \bar{z}) + q_2(z) \,(Q z \cdot \bar{z}) +q_3(z) \,(R z \cdot \bar{z})=0,\\
r_1(z)\,(P z \cdot \bar{z}) + r_2(z) \,(Q z \cdot \bar{z}) +r_3(z) \,(R z \cdot \bar{z})=0.
 \end{eqnarray}
Или же в векторной форме
\begin{eqnarray} \label{dv3}
\nonumber
 p_1(z)\,P z  + p_2(z) \,Q z  +p_3(z) \,R z =0,\\
\nonumber
 q_1(z)\,P z  + q_2(z) \,Q z  +q_3(z) \,R z =0,\\
 r_1(z)\,P z  + r_2(z) \,Q z  +r_3(z) \,R z =0.
 \end{eqnarray}
 Запишем первую координату первого соотношения (\ref{dv3}), вторую  -- второго, третью -- третьего,   получаем
 \begin{eqnarray*}
p_1^2+p_2 \, q_1+p_3 \, r_1=0,\\
q_1 \, p_2 + q_2^2+q_3 \, r_2=0,\\
r_1 \, p_3 +r_2 \, q_3+r_3^2=0.
 \end{eqnarray*}
 Откуда с помощью леммы 5 получаем
 \begin{eqnarray*}
 \nonumber
 p_1(z)=\alpha_2 \, p_2(z)+\alpha_3 \, p_3(z), \\
 q_2(z)=\beta_1 \, q_1(z)+\beta_3 \, q_3(z), \\
 r_3(z)=\gamma_1 \, r_1(z)+\gamma_2 \, r_2(z).
 \end{eqnarray*}
Тогда  (\ref{dv3}) принимает вид
\begin{eqnarray} \label{dv31}
\nonumber
 p_2(z)\,(Q+\alpha_2\,P) \,z   + p_3(z)\,(R+\alpha_3\,P) \,z =0,\\
\nonumber
q_1(z)\,(P+\beta_1\,Q) \,z   + q_3(z)\,(R+\beta_3\,Q) \,z =0,\\
r_1(z)\,(P+\gamma_1\,R) \,z   + r_2(z)\,(Q+\gamma_2\,R) \,z =0
  \end{eqnarray}
Из линейной независимости $(P,Q,R)$ следует, что в каждой паре
 $$(p_2,p_3),  (q_1,p_3),(r_1,r_2)$$
 либо обе формы одновременно равны нулю, либо непропорциональны.

 {\it Случай 3.0.}  Все три пары обратиться в ноль не могут, т.к. это означало бы что все три первых базисных вектора попали в общее ядро. Таким, образом одна из трех пар - ненулевая и непропорциональная. Пусть это $(r_1,r_2)$.

 {\it Случай 3.1.} Пусть обе оставшиеся пары -- нулевые.  Запишем первую координату третьего соотношения  (\ref{dv31}), получим $\gamma_1 \, r_1^2 +\gamma_2 \, r_1 \, r_2=0$. Откуда следует, что $\gamma_1=\gamma_2=0$. Тогда с помощью леммы 3 и невырожденности получаем, что $P\,z=r_2(z)\,A, \; Q\,z=-r_1(z)\,A, \; A \neq 0$, а из леммы 6, что как $r_1$, так и $r_2$ пропорциональны $ ( z\cdot \bar{A})$. Противоречие.

{\it Случай 3.2} Пусть только одна пара обратилась в ноль, скажем $p_1=p_2=0$, а формы $(q_1,q_3)$, также как и $(r_1,r_2)$ -- попарно непропорциональны. Запишем первую координату второго соотношения (\ref{dv31}), получим
$\beta_1 \,q_1^2+q_3\,r_1+\beta_3 \, q_3 \, q_1=0.$   Откуда сразу следует, что $r_1=\nu \,q_1$.
Запишем первую координату третьего соотношения (\ref{dv31}), получим
$\gamma_1 \,r_1^2+r_2\,(q_1+\gamma_2 \,r_1)=0.$   Откуда, ввиду того что $r_2$ непропорциональна $r_1$ следует, что $\gamma_1=0,\;
\gamma_2=-1$, т.е. $r_3=-r_2$. Из соотношения $r_1(z) \, P \,z+r_2(z)\,(R-Q)\,z=0$ следует, что
$r_j=q_j, \; j=1,2,3, \,\beta_1=0,\; \beta_3=-1,$ т.е. $q_2=-q_3$. Откуда получаем, что $r_1(z) \, P \,z+r_2(z)\,(Q-R)\,z=0$, т.е. $R=Q$. Противоречие.

\vspace{3ex}

{\it Случай 3.3.} Пусть все три пары попарно непропорциональны. Тогда из (\ref{dv31}), леммы 3 и невырожденности получаем
\begin{eqnarray} \label{dv32}
\nonumber
(Q+\alpha_2\,P)\,z = -p_3(z)\,A, \quad (R+\alpha_3\,P)\,z = p_2(z)\,A,\\
\nonumber (P+\beta_1\,Q)\,z = -q_3(z)\,B, \quad (R+\beta_3\,Q)\,z = q_1(z)\,B,\\
(P+\gamma_1\,R)\,z = -r_2(z)\,C, \quad (Q+\gamma_2\,R)\,z = r_1(z)\,C,
\end{eqnarray}
где $A,B,C$ -- три ненулевых вектора.
Выписывая для первого равенства третью координату, для второго -- вторую, для третьего -- третью, для четвертого -- первую,
для пятого -- вторую и для шестого -- первую, получим следующий набор соотношений
\begin{eqnarray}\label{eq1}
q_3(z)=-(\alpha_2+a_3)\,p_3(z),     \quad (\alpha_2+a_3)(\beta_1+b_3)=1,\\
\nonumber r_2(z)=(a_2-\alpha_3)\,p_2(z),     \quad (\alpha_3-a_2)(\gamma_1+c_2)=1,\\
\nonumber r_1(z)=(b_1-\beta_3)\,q_1(z),     \quad (b_1-\beta_3)(c_1-\gamma_2)=1.
\end{eqnarray}

Таким образом все девять линейных форм $(p_1,p_2,p_3,q_1,q_2,q_3,r_1,r_2,r_3)$ являются линейными комбинациями трех $(q_1,p_2,p_3)$.
\begin{eqnarray*}
p_1=\alpha_2 \, p_2+\alpha_3 \, p_3,  \;\; q_2=\beta_1\, q_1-\beta_3\, (a_3+\alpha_2) p_3, \;\; q_3=-(\alpha_2+a_3)\,p_3,\\
r_1=(b_1-\beta_3)\,q_1, \;\; r_2=(a_2-\alpha_3)\,p_2, \;\; r_3=\gamma_1\,(b_1-\beta_3)\,q_1+\gamma_2\,(a_2-\alpha_3)\,p_2.
\end{eqnarray*}
А также получаем, что
$$\beta_1 = \frac{1}{\alpha_2+a_3}-b_3, \;\;  \gamma_1 = \frac{1}{\alpha_3-a_2}-c_2, \;\;
\gamma_2= \frac{1}{\beta_3-b_1}+c_1,$$
где знаменатели не обращаются в ноль.

В таком случае все восемнадцать соотношений, полученных из первых трех координат (\ref{dv32}) представляют собой
систему $L(q_1,p_2,p_3)=0$ линейных соотношений между $(q_1,p_2,p_3)$.
Коэффициенты этой системы рационально зависят от $(a_1,a_2,a_3,b_1,b_2,b_3,c_1,c_2,c_3)$ и $(\alpha_2,  \alpha_3,\beta_3)$,
знаменатели в силу (\ref{eq1}) в ноль не обращаются.

Пронумеруем эти уравнения подряд, т.е.  первая координата первой группы, вторая координата первой группы, третья координата первой группы, первая координата второй группы и так до третьей координаты шестой группы. Получим
$$ l^1_j \, q_1+  l^2_j \, p_2+l^3_j \,p_3=0, \quad j=1,\dots,18.$$
Учитывая, что ранг пары форм $(p_2,p_3)$ равен двум, ранг полученной матрицы коэффициентов размером $3 \times18$ равен единице. Получаем систему из 34-х уравнений
\begin{eqnarray} \label{34}
\nonumber
e_{2j-3}=l^2_j \, l^1_1-l^1_j \, l^2_1, \quad e_{2j-2}=l^3_j \, l^1_1-l^1_j \, l^3_1, \quad  \quad j=2,\dots,18.
\end{eqnarray}
Пусть $E_m$ - числитель выражения $e_m$. Получаем систему из 34-х полиномиавльных соотношений на 12 переменных.
$E_m=0, \, m=1, \dots,34$.  Часть соотношений (\ref{dv32}) мы использовали выше. Непосредственно убеждаемся, что
$$E_3= E_4=E_7=E_8=E_{15}=E_{17}=E_{18}=E_{26}=0.$$
Систему из оставшихся 26 уравнений обозначим
\begin{equation}\label{L}
\mathcal{L}(a_1,a_2,a_3,b_1,b_2,b_3,c_1,c_2,c_3,\alpha_2,  \alpha_3,\beta_3)=0.
\end{equation}

\vspace{3ex}

{\bf Лемма 18:}  (a) Решения системы (\ref{L}) это объединение трех неприводимых компонент
\begin{eqnarray*}
\nonumber
 a_1=a_2=a_3=\alpha_2=\alpha_3=0,\\
 a_1=a_2=a_3=b_1=b_2=b_3=\alpha_3=0,\\
 a_1=a_2=a_3=b_1=b_2=b_3=c_1=c_2=c_3=\alpha_2 \, \beta_3+\alpha_3=0. \\
\end{eqnarray*}
 (b) Ни одно из этих решений не совместимо с условием случая (3.3.)\\
 {\it Доказательство:}  Пункт (a) проверяется непосредственным анализом системы, который упрощается в связи с тем, что
 что среди уравнений имеются уравнения $a_3=a_2$ и $(c_2-c_1)(a_2-\alpha_3)=0$. Поскольку из  (\ref{eq1}) следует, что
 $(a_2-\alpha_3) \neq 0$, то $c_2=c_1$.  Также это легко проверяется с использованием системы Maple.

 Далее, если решение принадлежит первой или второй компоненте, то $\alpha_3-a_2=0$, что противоречит (\ref{eq1}).
 Пусть решение из третьей компоненты. Тогда
 $$\beta_1=\frac{1}{\alpha_2}, \; \beta_3=-\frac{\alpha_3}{\alpha_2}, \;  \gamma_1 =\frac{1}{\alpha_3}, \; \gamma_2= -\frac{\alpha_2}{\alpha_3}.$$
При этом зависимые формы выражаются через $(q_1,p_2,p_3)$ следующим образом
$$q_2=-\alpha_2\,p_2, \; q_3=-\alpha_2\,p_3, \; r_1=-\alpha_2\, \alpha_3 \,p_2, -\alpha_3^2\,p_3, \; r_2=-\alpha_3\,p_2, \; r_3=-\alpha_3\,p_3.$$
Используем первые два равенства (\ref{dv32}), чтобы выразить $Q$ и $R$.
\begin{eqnarray}
\nonumber Q\,z=-\alpha_2 \,P\,z - p_3(z)\,A,\\
 R\,z=-\alpha_3 \,P\,z + p_2(z)\,A,
\end{eqnarray}
Результат подставим в оставшиеся четыре, получим.
\begin{eqnarray}\label{eq2}
\nonumber
B=-\frac{1}{\alpha_2}\,A, \; q_1=-\alpha_2^2\,p_2-\alpha_2\,\alpha_3 \,p_3, \; C=\frac{1}{\alpha_3^2}\,A, \; \frac{\alpha_2}{\alpha_3}\,p_2=0.
\end{eqnarray}
Последнее равенство не может быть выполнено т.к. $\alpha_2 \neq 0$. Это противоречие завершает наше рассуждение.

\vspace{3ex}

{\bf Теорема 19:}  Если $\mathcal{Q}$ -- невырожденная модельная квадрика коразмерности $k=3$,  то разложение
${\rm aut}\,\mathcal{Q}$ на весовые  компоненты имеет вид
$$g_{-2}+g_{-1}+g_0+g_1+g_2.$$
{\it Доказательство:} Как хорошо известно, подалгебра $g_{-}=g_{-2}+g_{-1}$
для любой невырожденной квадрики -- фундаментальна по Танаке (\cite{SS18}, \cite{G18}).  Выше нами было показано, что $g_3=0$. Откуда сразу следует тривиальность и всех компонент с б$\acute{o}$льшими весами. Теорема доказана.\\

\vspace{3ex}

Рассмотренные выше доказательства неисключительности квадрик малой коразмерности использовали только второе соотношение
(\ref{88a}).  Квадрикам малых коразмерностей противостоят, в некотором смысле, квадрики очень высоких коразмерностей. Самая большая коразмерность при фиксированном $n$ -- это $k=n^2$. Это единственная квадрика, которая задается базисным набором в пространстве эрмитовых форм. Введем следующие обозначения для переменных группы $w$.
\begin{eqnarray*}
w_{\alpha\alpha}=u_{\alpha\alpha}+i\,v_{\alpha\alpha}, \quad 1 \leq   \alpha \leq n.\\
w^R_{\alpha\beta}=u^R_{\alpha\beta}+i\,v^R_{\alpha\beta}, \quad w^I_{\alpha\beta}=u^I_{\alpha\beta}+i\,v^I_{\alpha\beta}, \quad 1 \leq \beta < \alpha \leq n.\\
\end{eqnarray*}
Тогда уравнения квадрики примут вид
\begin{eqnarray*}
v_{\alpha\alpha}=z_\alpha \, \bar{z}_\alpha, \quad 1 \leq   \alpha \leq n.\\
v^R_{\alpha\beta}=2\,{\rm Re} \, z_\alpha \, \bar{z}_\beta,  \quad  v^I_{\alpha\beta}=2\,{\rm Im} \, z_\alpha \, \bar{z}_\beta, \quad
1 \leq \beta < \alpha \leq n.
\end{eqnarray*}

\vspace{3ex}
{\bf Утверждение 20:} Если $\mathcal{Q}$ -- невырожденная квадрика коразмерности $k=n^2$ (последняя квадрика), то $g_{+}=0$.  Другими словами такая квадрика является жесткой.\\
{\it Доказательство:} Будем выделять некоторые координаты в первом векторном соотношении (\ref{88a}). Запишем $(\alpha\alpha)$-координату
$$ A_\alpha(u)(z,z) \, \bar{z}_\alpha = 2\, i\, z_\alpha\, \Delta \bar{a}_\alpha(u).$$
Из делимости на $\bar{z}_\alpha$ следует, что $a_\alpha$ зависит только от $u_{\alpha\alpha}$, причем
$A_\alpha(u)(z,z)= 2\, i\, z^2_\alpha\, \bar{a}'_\alpha(u_{\alpha\alpha})$. Запишем $(\alpha 1)^R$-координату (\ref{88a})  для $\alpha >1$, получим
$$ \bar{a}'_\alpha(u_{\alpha\alpha}) \, z^2_\alpha\,\bar{z}_1+  \bar{a}'_1(u_{11}) \, z^2_1\,\bar{z}_\alpha=
 z_\alpha\, \bar{a}'_1 |z_1|^2+z_1 \, \bar{a}'_\alpha \, |z_\alpha|^2$$
Получаем, что $a=const, \; A=0.$   Откуда сразу следует, что $g_1=0$, а в силу фундаментальности, что и $g_{+}=0$.   Утверждение доказано.

\vspace{3ex}

Отметим, что здесь, в противоположность малым коразмерностям, мы использовали только  первое соотношение  (\ref{88a}), причем небольшую его часть,  и совсем не использовали второго соотношения (\ref{88a}). Это наблюдение позволяет предположить,
что для больших коразмерностей основные ограничения на $(a,A)$ следуют из первого соотношения, а для малых -- из второго.

\vspace{3ex}

{\bf 4.  $RAQ$-квадрики}

\vspace{3ex}

В работе \cite{ES97} был рассмотрен  весьма интересный класс квадрик, т.ч. $n=k$ (коразмерность равна $CR$-размерности) -- $RAQ$-квадрики. Такие квадрики находятся во взаимно-однозначном соответсявии с конечномерными вещественными коммутативными (ассоциативными)  алгебрами с единицей .

Пусть $\mathcal{A}$ -- такая алгебра размерности $n$.  Если $X$ и $Y$ -- элементы $\mathcal{A}$, то через $x \cdot y$ обозначим их произведение как элементов алгебры.  Пусть  $\mathcal{A}^{c}=\mathcal{A} \bigotimes \mathbf{C}$ -- комплексификация $\mathcal{A}$. Если $Z \in  \mathcal{A}^{c}$,   то через $\bar{Z}$ обозначим комплексное сопряжение, канонически определенное в $\mathcal{A}^{c}$.  Квадрика, соответствующая алгебре $\mathcal{A}$ имеет вид

\begin{equation}\label{RAQ}
Q=\{ (Z,W) \in (\mathcal{A}^{c})^2: \; {\rm Im} \, W = Z \cdot \bar{Z}\}.
\end{equation}

Если задать вопрос так:  когда квадрика типа $(n.n)$ имеет такой вид? То ответ -- это два условия на форму $<z,\bar{z}>$: \\
--  вещественность -- в некоторых координатах  форма на вещественных векторах принимает вещественные значения,\\
--  ассоциативность --  $<<p,q>,r>=<p,<q,r>>$ для всех  $p,\, q, \,r \in \mathbf{R}^n$.\\
Поэтому они и называются $RAQ$-квадриками - Real Associative Quadrics.

В \cite{ES97} было показано, что условие невырожденности для $RAQ$-квадрики равносильно наличию в $\mathcal{A}$ единицы.

\vspace{3ex}

{\bf Теорема 21:} Исключительных  $RAQ$-квадрик не существует.\\
{\it Доказательство:}
Для доказательства воспользуемся достаточным условием неисключительности - теоремой 6. Проверим выполнение  условий (I) и (II) .

Условие (I).  Пусть $(E_1,\dots,E_n)$ -- базисные элементы  алгебры $\mathcal{A}$, причем $E_1$ -- это единица алгебры. Тогда координатные операторы имеют вид $H_j \, Z=E_j \cdot Z$ и, тем самым,  $H_j \, E_1=E_j$,  а набор $(E_1,\dots,E_n)$ имеет ранг $n$.

Условие (II).   Образ отображения $(Z',Z) \rightarrow <Z',\bar{Z}>=Z' \cdot \bar{Z}$  совпадает с  $\mathcal{A}^c$, т.к. $(Z',E_1)$ переходит в $Z'$.  Теорема доказана.

\vspace{3ex}

{\bf Следствие  22:} Если $\mathcal{Q}$ -- это невырожденная $RAQ$-квадрика, то\\
(a)   ${\rm aut}\,  \mathcal{Q} = g_{-2}+g_{-1}+g_{0}+g_{1}+g_{2}$\\
(b)   Подгруппа $G_{+}$ -- подгруппа нелинейных автоморфизмов  $\mathcal{Q}$, сохраняющих начало координат, описывается формулой Пуанкаре (см.  \cite{ES97}), а именно
\begin{eqnarray}
\nonumber
Z^{*} = (Z + a \cdot W) \cdot (1 - 2 \,i \,\bar{a} \cdot  Z - (r + i \, a \cdot \bar{a} )\cdot W)^{-1},\\
W^{*}  = W \cdot (1 - 2 \,i\, \bar{a} \cdot Z - (r + i \,a \cdot \bar{a})\cdot W)^{-1},
\end{eqnarray}
где $a \in \mathcal{A}^c,  \; r \in \mathcal{A}$.

\vspace{5ex}

{\bf 5. Оценка на длину подалгебры}  $\mathbf{g_{+}}$

\vspace{3ex}
В \cite{BER} приведена теорема (теорема 12.3.11), позволяющая оценить номер струи $j$, от которой зависит отображение аналитического   $l$-невырожденного ростка коразмерности $k$, а именно
$$j \leq  (k+1) \, l.$$
Поскольку невырожденная квадрика $\mathcal{Q}$ является $1$-невырожденной и имеет кончный тип в каждой точке (в частности,  минимальна), то для степени $d$ коэффициентов полей из ${\rm aut}\,\mathcal{Q}$ получаем следующую оценку
$$d \leq  (k+1).$$
В этой части работы мы получим эту оценку и другие близкие результаты, пользуясь методами, отличными от методов работы \cite{BER} (множества Сегре). А именно, мы используем технику преобразования Фурье в пространстве обобщенных функций, которая лежит в основе доказательства теоремы Эренпрайса-Паламодова \cite{P}. Эта теорема была использована в работе \cite{VB88} при получении критерия конечномерности алгебры квадрики $\mathcal{Q}$.

\vspace{3ex}

 Пусть $F(s)=(F_1(s),\dots,F_n(s))$ -- обобщенная вектор-функция, которая является преобразованием Фурье  вектор-функции $a(u)=(a_1(u),\dots,a_n(u))$. Тогда, применяя преобразование Фурье ко второму уравнению (\ref{88a}), получаем, что для каждого $z \in \mathbf{C}^n$ выполнено соотношение
\begin{equation}\label{F1}
(s_1 \, <z,\bar{z}>_1 + \dots +s_k \, <z,\bar{z}>_k)^2 \, <F(s),\bar{z}>=0
\end{equation}
Пусть $C[s]=C[s_1,\dots,s_k]$ -- кольцо полиномов от $k$-мерной  переменной $s$ с комплексными коэффициентами, $C[u]=C[u_1,\dots,u_k]$ -- такое же кольцо полиномов от $k$-мерной  переменной $u$. Пусть $(C[s])^k$ -- свободный
модуль размерности $k$ над $C[s]$, а $(C[u])^n$ -- свободный модуль размерности $n$ над $C[u]$. Пусть
$$\{e_j=(0, \dots,0,1 \mbox{ (на j-м месте)},0, \dots,0), \; j=1,\dots,k \}$$
 -- набор образующих модуля $(C[s])^k$.

Введем в рассмотрение модуль  $M=M(\mathcal{Q})$ -- подмодуль модуля $(C[s])^k$,   порожденный всеми образующими вида
\begin{equation}\label{M}
(s_1 \, <z,\bar{z}>_1 + \dots +s_k \, <z,\bar{z}>_k)^2 \, (<\varphi(s),\bar{z}>_1 \, e_1+ \dots+<\varphi(s),\bar{z}>_k \, e_k),
\end{equation}
где $z \in \mathbf{C}^n, \; \varphi(s)=(\varphi_1(s),\dots, \varphi_n(s)) \in (C[s])^n$.
Будем называть $M(\mathcal{Q})$ {\it характеристическим} подмодулем квадрики $\mathcal{Q}$ .

\vspace{3ex}
Приведем набор утверждений, которые для уравнения
\begin{equation}\label{a}
  \Delta^2 \, <a(u), \bar{z}>=0
\end{equation}
непосредственно следуют из общей схемы \cite{P} (см. также \cite{STF} глава 10).

\vspace{3ex}

{\bf Утверждение:} Если квадрика $\mathcal{Q}$  невырождена, то\\
(a) Пространство решений (\ref{a}) -- это линейное подпространство  $\mathcal{L}_1$ пространства полиномов степени не выше чем $p <\infty$.\\
(b) Пространство решений (\ref{F1}) -- это линейное подпространство  $\mathcal{L}_2$ пространства, порожденного $\delta_0$ (дельта-функция с носителем в нуле) и ее производными порядка не выше чем $p$.\\
(c) Пространства $\mathcal{L}_1$ и $\mathcal{L}_2$ -- изоморфны, изоморфизм устанавливается преобразованием Фурье.\\
(d)  Характеристический подмодуль  имеет конечную коразмерность, т.е. если $M'=((C[s])^k /M )$, то  $\dim \,M'  $ как линейного пространства конечна. Причем $M'$ как линейное пространство естественно отождествляется с  пространством $\mathcal{L}_3$, которое получено заменой производной дельта-функции на соответствующий моном от переменных $s$.   Пространство $\mathcal{L}_3$ изоморфно $\mathcal{L}_2$ и отождествляется с фактор-пространством $M'$
в силу того, что любой элемент  $(C[s])^k$ однозначно представим в виде суммы элементов $\mathcal{L}_3$ и $M$.

\vspace{3ex}

Это утверждение показывает, что все характеристики  пространства $\mathcal{L}_1$ решений  (\ref{a})  (размерность, степень, ...)  совпадают   с аналогичными характеристиками пространства $\mathcal{L}_3$, которое предсталяет собой прямое дополнение $M$ до $(C[s])^k$.

Поскольку кольцо полиномов $C[s]$ -- нётерово, то модуль $M$ -- конечно порожден. В соответствии с (\ref{M}) это означает, что
у $M$ имеется конечная система образующих вида
$$\{(l_{\nu}(s))^2 \, (\lambda_{\nu}^1 \, e_1 + \dots + \lambda_{\nu}^k \, e_k)\},$$
где $l_{\nu}(s)$ -- линейные формы, а $\lambda \in \mathbf{C}^k$.

Пусть $M_j$ -- проекция $M$ на $j$-ю координату $(C[s])^k$.  Это идеал в  $C[s]$,  порожденный формами вида
$\{(l_{\nu}(s))^2 \, \lambda_{\nu}^j \}$, где $ \lambda_{\nu}^j \neq 0$.  Пусть $r_j$ -- это ранг соответствующего набора
линейных форм $\{l_{\nu}(s)\}$.

\vspace{3ex}

{\bf Утверждение  23:} (a) Для всех $j$ ранг $r_j=k$.\\
(b) Дополнение к  $M$ -- пространство $\mathcal{L}_3$  -- это подпространство пространства наборов полиномов от $s$, чья суммарная степень не превосходит$k$.  Пространство $\mathcal{L}_2$ -- это подпространство пространства линейных комбинаций дельта-функции и ее производных порядка не выше $k$. Пространство $\mathcal{L}_1$ -- это подпространство пространства наборов полиномов от $u$, чья суммарная степень не превосходит $k$.\\
(с)
{\it Доказательство:} Пусть найдется $j$,  т.ч. $r_j < k$.  Выберем из набора $\{l_{\nu}(s)\}$ линейно независимый набор максимального ранга
$(L_1(s), \dots, L_{k'}(s))$   $k' <k$. Сделаем в пространстве переменного $s$ невырожденную линейную замену  $s \rightarrow \tilde{s}$, т.ч.
$\tilde{s}_j=L_j(s)$ для $j=1, \dots,k'$. В новых переменных образующие $M_j$ -- это $(\tilde{s}_1)^2,\dots, (\tilde{s}_{k'})^2$. Дополнение
$M_j$ до $C[\tilde{s}]$ бесконечномерно. Действительно, оно содержит все бесконечномерное кольцо полиномов от последней переменной $C[\tilde{s}_k]$. А из этого сразу следует, что бесконечномерно дополнение  к $M$. Это противоречие доказывает пункт (a).\\
Итак, среди образующих идеала есть квадраты всех координат. Поэтому дополнение содержит лишь те полиномы, чья степень по любой переменной не превосходит единицы. Это, в частности, означает, что суммарная степень не превосходит $k$. Возвращаясь к старым переменным, мы сохраняем эту оценку. Итак, проекция дополнения на каждую координату не содержит полиномов степени выше $k$.
Оставшиеся утверждения -- это следствие описанных изоморфизмов. Утверждение доказано.

\vspace{3ex}

Из леммы 4 сразу следует что если $A(u)(z,z)$ -- решение (\ref{88a}),  то
$${\rm deg}_u \,A(u)(z,z) \leq k-1.$$
  Рассуждения, аналогичные тем, что были приведены выше, показывают, что если $(C(u),b(u))$ --  решение (\ref{88b}), то
  $${\rm deg}_u \,C(u) \leq k,   \quad    {\rm deg}_u \,b(u) \leq 2 \,k.$$
   Из этих оценок
можно получить общие оценки на степени коэффициентов произвольного элемента алгебры ${\rm aut}\,\mathcal{Q}$. В результате
для нечетной компоненты получаем оценку степени $(k+1)$, а для четной -- оценку $2 k$. Однако можно поступить иначе. Это даст более точную информацию, а к тому же не использует других оценок, кроме полученной выше оценки степени $a$.  Можно сначала дать оценку {\it веса} нечетной компоненты, а затем воспользоваться фундаментальностью алгебры.

\vspace{3ex}

{\bf Теорема 24:}  Если $\mathcal{Q}$ -- невырожденная модельная квадрика коразмерности $k$,  то \\
(a)  Вес полей нечетного веса из   ${\rm aut}\,\mathcal{Q}$ не превосходит $2\,k-1$.\\
(b)  Вес произвольного поля не превосходит $2\,k$, т.е.
$${\rm aut}\,\mathcal{Q}=g_{-2}+g_{-1}+g_{0}+g_{1}+\dots+g_{2k}$$
(c) Степени коэффициентов полей из ${\rm aut}\,\mathcal{Q}$ не превосходят $(k+1)$.\\
{\it Доказательство:} В силу оценок, полученных для $a$ и $A$, сразу получаем, что вес нечетной компоненты
\begin{eqnarray*}
2 \, {\rm Re} \, \left( (a(w)+A(w)(z,z)) \,\frac{\partial}{\partial z}
+ 2\,i\,<z,\bar{a}(w)> \,\frac{\partial}{\partial w}  \right)
\end{eqnarray*}
не превосходит $2k-1$. Это пункт (a). Из фундаментальности следует, что не существует ненулевых полей веса больше, чем $2k$. Это пункт (b). Как было отмечено, степени коэффициентов для полей нечетного веса не превосходят $(k+1)$.  Оценку для степени коэффициентов для полей четного веса получим из оценки веса. Если вес поля четного веса
\begin{eqnarray*}
2 \, {\rm Re} \, \left( C(w)\,z \,\frac{\partial}{\partial z}
+ b(w) \,\frac{\partial}{\partial w}  \right)
\end{eqnarray*}
не превосходит $2k$, то степени как $C(u)z$, так и  $b(u)$ не превосходит также  $(k+1)$.
Теорема доказана.

\vspace{3ex}

Как известно \cite{T88},  группа локальных автоморфизмов невырожденной квадрики -- это подгруппа группы бирациональных автоморфизмов объемлющего пространства $\mathbf{C}^{n+k}$ с равномерной оценкой на степень.   Полученная в предыдущей теореме оценка степени инфинетизимальных автоморфизмов позволяет дать оценку на степень автоморфизмов.

\vspace{3ex}

{\bf Следствие  25:}  Пусть $\mathcal{Q}$ -- невырожденная квадрика типа $(n,k)$, а ${\rm Aut} \, \mathcal{Q}$ -- локальная группа
ее  автоморфизмов, голоморфных в окрестности нуля. Тогда  ${\rm Aut} \, \mathcal{Q}$ состоит из бирациональных преобразований
$\mathbf{C}^{n+k}$ чья степень не превосходит $(3n +3k+2)(k+1)$.

\vspace{3ex}

{\bf 6. Открытые вопросы}

\vspace{3ex}

  Если $(n,k)$ -- $CR$-тип ($n$ -- $CR$-размерность, $k$ -- коразмерность), то невырожденные квадрики возможны в диапазоне $1 \leq k \leq n^2$.\\

  {\bf Вопрос  26.}  Для каких $CR$-типов из этого диапазона существуют исключительные квадрики? \\

Как мы видели, кроме упомянутых случаев $k=1,2,3$, в которых исключительных квадрик  нет, к числу таких коразмерностей можно добавить $k=n^2$   (утверждение 20).\\

 Выше (теорема 7) было показано, что не существует квадрик коразмерности 4 для $CR$-размерности $n \leq 3$.  Пример исключительной квадрики коразмерности $k=4$, приведенный  в \cite{GM}, имеет  $CR$-размерность $n =6$.  В качестве частного случая вопроса 26 можно
 предложить следующий вопрос.  Какова минимальная $CR$-размерность допускающая исключительные квадрики коразмерности 4?  (Варианты ответа: $4,\, 5, \, 6$.)

\vspace{2ex}

 Выше (теорема 24) было показано, что веса положительных компонент алгебры Ли невырожденной квадрики  коразмерности $k$ не превосходят $2k$. Но эта оценка не подкреплена примером и вопрос о точной оценке остается открыт.\\

 {\bf Вопрос  27:}  Существуют ли невырожденные квадрики, т.ч. длина $g_{+}$ равна $2k$?\\

\vspace{2ex}

 Все известные нам примеры не жестких квадрик это примеры, где старший вес -- четный. \\

{\bf Вопрос  28:} Существуют ли невырожденные квадрики, т.ч. длина $g_{+}$  равна нечетному числу?\\

\vspace{2ex}

     В тех $(n,k)$-случаях, когда исключительные квадрики появляются, что предсталяет собой подмножество исключительных квадрик в совокупности всех невырожденных квадрик данного типа? Как мы отмечали выше исключительные квадрики задаются конечной совокупностью полиномиальных соотношений типа "равно нулю" и "не равно нулю" (полуалгебраическое множество). При этом естественная точка зрения состоит в том, чтобы проводить это исследование не в пространстве $\mathcal{H}$  наборов $k$ эрмитовых форм, а в пространстве модулей \cite{VB06}, т.е. после факторизации этого пространства по известному линейному действию. Интерес представляет ответ на следующий вопрос.

 {\bf Вопрос  29:}  Какова коразмерность подпространства модулей исключительных квадрик в пространстве модулей всех невырожденных квадрик?\\

\end{document}